\documentclass[a4,12pt]{article}%
\usepackage{amsmath}
\usepackage{graphicx}
\usepackage{hyperref}
\usepackage{eurosym}
\usepackage{amsfonts}
\usepackage{amssymb}%
\providecommand{\U}[1]{\protect\rule{.1in}{.1in}}
\providecommand{\U}[1]{\protect\rule{.1in}{.1in}}

\begin{document}

\begin{center}
{\huge A Multidimensional Fatou Lemma for }\bigskip

{\huge Conditional Expectations }

\bigskip\ \bigskip

E. Babaei$^{\ast}$, I. V. Evstigneev$^{\ast\ast}$ and K. R.
Schenk-Hopp\'{e}$^{\ast\ast\ast}$\bigskip
\end{center}

{\small \noindent\textbf{Abstract:} The classical multidimensional version of
Fatou's lemma (Schmeidler \cite{Schmeidler1970}) originally obtained for
unconditional expectations and the standard non-negative cone in a
finite-dimensional linear space\ is extended to conditional expectations and
general closed pointed cones.\bigskip}

{\small \noindent\textbf{Key words and Phrases:} Cones in linear spaces;
Induced partial orderings; Sequences of random vectors, Fatou's lemma;
Conditional expectations\bigskip}


{\small \bigskip\bigskip\bigskip\bigskip\bigskip\bigskip\bigskip
\bigskip\bigskip\bigskip\bigskip\bigskip\bigskip\bigskip}

\rule{5cm}{0.01cm}

$^{\ast}${\small Department of Economics, University of Manchester, Oxford
Road, Manchester M13 9PL, UK. E-mail:
esmaeil.babaeikhezerloo@manchester.ac.uk.}

$^{\ast\ast}${\small Department of Economics, University of Manchester, Oxford
Road, Manchester M13 9PL, UK. E-mail: igor.evstigneev@manchester.ac.uk.
(Corresponding author.)}

$^{\ast\ast\ast}${\small Department of Economics, University of Manchester,
Oxford Road, Manchester M13 9PL, UK, and Department of Finance, NHH --
Norwegian School of Economics, Helleveien 30, 5045 Bergen, Norway. E-mail:
klaus.schenk-hoppe@manchester.ac.uk.}\pagebreak

\textbf{1.} Fatou's lemma in several dimensions, the first version of which
was obtained by Schmeidler \cite{Schmeidler1970}, is a powerful
measure-theoretic tool initially developed in Mathematical Economics in
connection with models of "large" economies with atomless measure spaces of
agents; see Aumann \cite{Aumann1966} and Hildenbrand \cite{Hildenbrand1974}.
In this note we provide two new versions of this lemma: one for unconditional
and the other for conditional expectations. Both deal with cones in an
$n$-dimensional linear space $\mathbb{R}^{n}$ more general than the
non-negative orthant $\mathbb{R}_{+}^{n}$ as considered in
\cite{Schmeidler1970}. Our results are motivated by the applications of the
theory of von Neumann-Gale \cite{vonNeumann1937,Gale1956} dynamical systems to
the modeling of financial markets with frictions---transaction costs and
portfolio constraints
\cite{DempsterEvstigneevTaksar2006,EvstigneevSchenkHoppe2008,EvstigneevZhitlukhin2013,BabaeiEvstigneevPirogov2018,BESH2018}%
.

\textbf{2.} Let $(\Omega,\mathcal{F},P)$ be a probability space and $C$ a
pointed closed cone\footnote{A set $C$ in a linear space is called a
\textit{cone} if it contains with any its elements $x,y$ any non-negative
linear combination $\lambda x+\mu y$ ($\lambda,\mu\geq0$) of these elements.
The cone $C$ is called \textit{pointed }if the inclusions $x\in C$ and $-x\in
C$ imply $x=0$.} in $\mathbb{R}^{n}$. We write $a\leq_{C}b$ if $b-a\in C$. Let
$|\cdot|$ be a norm in $\mathbb{R}^{n}$. The distance measured in terms of the
norm $|\cdot|$ between a point $a$ and a set $A$ in $\mathbb{R}^{n}$ is
denoted by $\rho(a,A)$. We will use the standard notation Ls$(x_{k})$ for the
set of limit points of the sequence $(x_{k})$.

Recall that a sequence of random variables $\beta_{k}(\omega)$ is called
\textit{uniformly integrable} if
\begin{equation}
\lim_{H\rightarrow\infty}\sup_{k}E|\beta_{k}|\mathbf{1}_{\{|\beta_{k}|\geq
H\}}=0. \label{unif-integr}%
\end{equation}
Property (\ref{unif-integr}) holds if and only if the following two conditions
are satisfied: a) $\sup E|\beta_{k}|<\infty$; b) $\lim E|\beta_{k}%
|\mathbf{1}_{\Gamma_{k}}=0$ for any sequence of events $\Gamma_{k}$ with
$P(\Gamma_{k})\rightarrow0$ (see, e.g., Neveu \cite{Neveu1965}).

\textbf{Theorem 1.} \textit{Let }$x_{k}(\omega)$\textit{,}$\;k=1,2,...$%
\textit{, be a sequence of random vectors in }$\mathbb{R}^{n}$ \textit{such
that }$E|x_{k}(\omega)|<\infty$\textit{\ and }$Ex_{k}(\omega)\rightarrow y
$\textit{, where }$y$\textit{\ is some vector in} $\mathbb{R}^{n}%
$.\textit{\ If the sequence }$\rho(x_{k}(\omega),C)$\textit{\ is uniformly
integrable, then there exist integer-valued random variables }%
\begin{equation}
1<k_{1}(\omega)<k_{2}(\omega)<... \label{k-m}%
\end{equation}
\textit{\ and a random vector }$x(\omega)$\textit{\ such that }$E|x|<\infty
$\textit{, }%
\begin{equation}
\lim_{m\rightarrow\infty}x_{k_{m}(\omega)}(\omega)=x(\omega
)\text{\ \textit{(a.s.)}} \label{Ls-conv}%
\end{equation}
\textit{\ and}%
\[
Ex(\omega)\leq_{C}y.
\]

Theorem 1 is a version of the multidimensional Fatou lemma in
\cite{Schmeidler1970} where it is assumed that $C=\mathbb{R}_{+}^{n}$ and
$x_{k}(\omega)\in C$, so that $\rho(x_{k}(\omega),C)=0$. It also extends a
result in the paper by Cornet et al. \cite{Cornet-et-al2003}, Theorem B,
p.\ 194, in which the function\ $\rho(x_{k}(\omega),C)$\ is required to be
integrably bounded.

\textbf{3. }\textit{Proof of Theorem 1.} \textit{1st step.} We have
$x_{k}(\omega)=c_{k}(\omega)+b_{k}(\omega)$, where $c_{k}(\omega)\in C$ and
$\left\vert b_{k}(\omega)\right\vert =\rho(x_{k}(\omega),C)$ (a.s.). By
assumption, the sequence $\left\vert b_{k}\right\vert $ is uniformly
integrable, and so $H:=\sup E\left\vert b_{k}\right\vert <\infty$. Since
$c_{k}=x_{k}-b_{k}$, we have $E|c_{k}|\leq E|x_{k}|+H<\infty$, and so the
random vectors $c_{k}$ are integrable. Furthermore, the sequence
$Ec_{k}=Ex_{k}-Eb_{k}$ is bounded because $\sup|Eb_{k}|\leq\sup E|b_{k}|=H$
and the sequence $Ex_{k}$ is bounded since it converges.

Note that the boundedness of $Ec_{k}$ implies the boundedness of $E|c_{k}|$
because the random vectors $c_{k}(\omega)$ take on their values in the closed
pointed cone $C$. Indeed, consider a strictly positive linear functional $g$
on $C$ ($gc>0$, $0\neq c\in C$); it exists for each closed pointed cone. For
such a functional $g$, there exists $\gamma>0$ such that $gc\geq\gamma|c|$ for
all $c\in C$. Consequently, $gEc_{k}=Egc_{k}\geq\gamma E|c_{k}|$, which proves
the boundedness of $E|c_{k}|$. This, in turn, implies that $E|x_{k}|$ is
bounded because $\sup E|x_{k}|\leq\sup(E|c_{k}|+E|b_{k}|)\leq\sup
E|c_{k}|+H<\infty$.

\textit{2nd step.} Since $\sup E|x_{k}|<\infty$, we can use the "biting lemma"
(e.g. Saadoune and Valadier \cite{SaadouneValadier1995}, p.\ 349) and find a
subsequence $(x_{k_{l}})$ of $(x_{k})$ and measurable sets $\Gamma
_{1}\supseteq\Gamma_{2}\supseteq...$ such that $%
{\textstyle\bigcap\nolimits_{l}}
\Gamma_{l}=\emptyset$ and the sequence $x_{k_{l}}^{\prime}:=x_{k_{l}%
}\mathbf{1}_{\Omega\backslash\Gamma_{l}}$ is uniformly integrable. Put
$x_{k_{l}}^{\prime\prime}:=x_{k_{l}}\mathbf{1}_{\Gamma_{l}}$. Clearly
$x_{k_{l}}=x_{k_{l}}^{\prime}+x_{k_{l}}^{\prime\prime}$.

For each $\omega$ the sequence $x_{k_{l}}^{\prime}(\omega)$ coincides with
$x_{k_{l}}(\omega)$ from some $l(\omega)$ on because every $\omega$ belongs to
$\Omega\backslash\Gamma_{l}$ from some $l(\omega)$ on. Consequently, for all
$\omega$ we have Ls$(x_{k_{l}}^{\prime})=$ Ls$(x_{k_{l}})\subseteq$
Ls$(x_{k})$.

The sequences $E|x_{k_{l}}^{\prime}|$ and $E|x_{k_{l}}^{\prime\prime}|$ are
bounded because $|x_{k_{l}}^{\prime}|\leq|x_{k_{l}}|$ and $|x_{k_{l}}%
^{\prime\prime}|\leq|x_{k_{l}}|$. By passing to a subsequence, we can assume
without loss of generality that\ $Ex_{k_{l}}^{\prime}\rightarrow y^{\prime}$
and $Ex_{k_{l}}^{\prime\prime}\rightarrow y^{\prime\prime}$ for some
$y^{\prime},y^{\prime\prime}\in\mathbb{R}^{n}$. Clearly, $y^{\prime}%
+y^{\prime\prime}=y$ because $Ex_{k_{l}}^{\prime}+Ex_{k_{l}}^{\prime\prime
}=Ex_{k_{l}}\rightarrow y$.

\textit{3rd step.} Since the sequence $(x_{k_{l}}^{\prime})$ is uniformly
integrable and $Ex_{k_{l}}^{\prime}\rightarrow y^{\prime}$, by Artstein's
theorem \cite{Artstein1979}, Theorem A, there exists an integrable random
vector $x(\omega)$ such that $x(\omega)\in$Ls$(x_{k_{l}}^{\prime}%
(\omega))\subseteq$Ls$(x_{k}(\omega))$ (a.s.) and $Ex(\omega)=y^{\prime}$. We
have%
\[
Ex_{k_{l}}^{\prime\prime}=Ex_{k_{l}}\mathbf{1}_{\Gamma_{l}}=Ec_{k_{l}%
}\mathbf{1}_{\Gamma_{l}}+Eb_{k_{l}}\mathbf{1}_{\Gamma_{l}}\rightarrow
y^{\prime\prime},
\]
where $Eb_{k_{l}}\mathbf{1}_{\Gamma_{l}}\rightarrow0$ because $P(\Gamma
_{l})\rightarrow0$ and the sequence $b_{k_{l}}$ is uniformly integrable. Thus
$Ec_{k_{l}}\mathbf{1}_{\Gamma_{l}}\rightarrow y^{\prime\prime}$. We have
$c_{k_{l}}(\omega)\mathbf{1}_{\Gamma_{l}}(\omega)\in C$ because $c_{k_{l}%
}(\omega)\in C$ and $0\in C$. Consequently, $Ec_{k_{l}}\mathbf{1}_{\Gamma_{l}%
}\in C$ as the set $C$ is convex (see, e.g., \cite{ArkinEvstigneev1987},
Appendix II, Lemma 1). Therefore $y^{\prime\prime}\in C$ since $Ec_{k_{l}%
}\mathbf{1}_{\Gamma_{l}}\rightarrow y^{\prime\prime}$ and $C$ is closed.

\textit{4th step.} We obtained that $y-y^{\prime}=y^{\prime\prime}\in C$,
i.e., $y^{\prime}\leq_{C}y$. Furthermore, $Ex(\omega)=y^{\prime}$, so that
$Ex(\omega)\leq_{C}y$, where $x(\omega)\in$Ls$(x_{k}(\omega))$ (a.s.). It
remains to observe that the inclusion $x(\omega)\in$Ls$(x_{k}(\omega))$ (a.s.)
implies the existence of a sequence $(k_{m}(\omega))_{m=1}^{\infty}$ of
integer-valued random variables such that (\ref{Ls-conv}) holds. Indeed, since
$x(\omega)\in$Ls$(x_{k}(\omega))$ (a.s.), for almost all $\omega$ there exists
a sequence $\kappa=(k_{m})_{k=1}^{\infty}$ of natural numbers $k_{m}$ for
which
\begin{equation}
1<k_{1}<k_{2}<...\ \text{and }\ \lim x_{k_{m}}=x(\omega). \label{x-k-m}%
\end{equation}
Denote by $A$ the set of $(\omega,\kappa)$ satisfying (\ref{x-k-m}). This set
is measurable with respect to $\mathcal{F}\times\mathcal{B}(\mathbb{N}%
^{\infty})$, where $\mathbb{N}^{\infty}:=\mathbb{N}\times\mathbb{N}\times...$
is the product of a countable number of copies of the discrete space
$\mathbb{N}:=\{1,2,...\}$ and $\mathcal{B}(\mathbb{\cdot})$ stands for the
Borel $\sigma$-algebra. Since $(\mathbb{N}^{\infty},\mathcal{B}(\mathbb{N}%
^{\infty}))$ is a standard measurable space\footnote{A measurable space is
called \textit{standard }if it is isomorphic to a Borel subset of a complete
separable metric space with the Borel measurable structure.}, we can apply
Aumann's measurable selection theorem (see e.g. \cite{ArkinEvstigneev1987},
Appendix I, Corollary 3) and construct a measurable mapping $\kappa(\omega)$
of $\Omega$ into $\mathbb{N}^{\infty}\ $for which $(\omega,\nu(\omega))\in A$
for almost all $\omega$. The sequence\ $\kappa(\omega)=(k_{m}(\omega
))_{m=1}^{\infty}$ of measurable integer-valued random variables satisfies
(\ref{x-k-m}) and (\ref{Ls-conv}). \hfill$\square$

\textbf{4. }Let $\mathcal{G}$ be a sub-$\sigma$-algebra of $\mathcal{F}$ and
let $C(\omega)$ be a pointed closed convex cone in $\mathbb{R}^{n}$ depending
$\mathcal{G}$-measurably\footnote{A set $C(\omega)\subseteq\mathbb{R}^{n}$ is
said to depend $\mathcal{G}$-measurably on $\omega$ if its graph
$\{(\omega,c):~c\in C(\omega)\}$ belongs to $\mathcal{G}\times\mathcal{B(}%
\mathbb{R}^{n})$.} on $\omega$. A random vector $x(\omega)$ is said to be
\textit{conditionally integrable} (with respect to the $\sigma$-algebra
$\mathcal{G}$) if $E[|x(\omega)|\,|\mathcal{G}]<\infty$ (a.s.). A sequence of
random variables $\beta_{k}(\omega)$, $k=1,2,...$, is said to be
\textit{uniformly conditionally integrable} if
\begin{equation}
\lim_{H\rightarrow\infty}\sup_{k}E[|\beta_{k}|\mathbf{1}_{\{|\beta_{k}|\geq
H\}}|\mathcal{G}]=0\ \text{(a.s.)}. \label{cond-unif-integr}%
\end{equation}

The following result is a version of Theorem 1 for conditional expectations.

\textbf{Theorem 2.} \textit{Let }$x_{k}(\omega)$\textit{,}$\;k=1,2,...$%
\textit{, be conditionally integrable random vectors in }$\mathbb{R}^{n}$
\textit{and }$y(\omega)$ \textit{a random vector in }$\mathbb{R}^{n}%
$\textit{\ such that }%
\begin{equation}
E[x_{k}(\omega)|\mathcal{G}]\rightarrow y(\omega)\mathit{\ }%
\text{\textit{(a.s.).}} \label{cond-conv}%
\end{equation}
\textit{\ If the sequence }$\rho(x_{k}(\omega),C(\omega))$\textit{\ is
uniformly conditionally integrable, then there exists a sequence of
integer-valued random variables }$1<k_{1}(\omega)<k_{2}(\omega)<...$%
\textit{\ and a conditionally integrable random vector }$x(\omega
)$\textit{\ such that }%
\[
\lim_{m\rightarrow\infty}x_{k_{m}(\omega)}(\omega)=x(\omega
)\text{\ \textit{(a.s.)}}%
\]
\textit{and}%
\[
E[x(\omega)|\mathcal{G}]\leq_{C(\omega)}y(\omega)\mathit{\;}%
\text{\textit{(a.s.).}}%
\]

In the case when $C(\omega)=\mathbb{R}_{+}^{n}$ Theorem 2 was proved in
\cite{BahsounEvstigneevTaksar2008}, Appendix A, Proposition A.2. For reviews
of various results related to multidimensional Fatou lemmas, see Balder and
Hess \cite{Balder-Hess1995} and Hess \cite{Hess2002}.

Some comments regarding the assumptions of Theorem 2 are in order. Clearly a
sequence of random variables $\beta_{k}(\omega)$ is uniformly conditionally
integrable if it is \textit{conditionally integrably bounded,} i.e.
$|\beta_{k}(\omega)|\leq\alpha(\omega)$, where $E[|\alpha(\omega
)||\mathcal{G]<\infty}$ (a.s.). The last condition holds, in particular, if
$\beta_{k}(\omega)$ is (unconditionally) integrably bounded: $|\beta
_{k}(\omega)|\leq\alpha(\omega)$ (a.s.) where $E|\alpha(\omega
)|\mathcal{<\infty}$. It should be noted that uniform integrability does not
necessarily imply uniform conditional integrability.

\textbf{5.} \textit{Proof of Theorem 2.} \textit{1st step.} Let us regard the
sequence of random vectors $x^{\infty}(\omega):=(x_{1}(\omega),x_{2}%
(\omega),...)$ as a random element of the standard measurable space
$(X^{\infty},\mathcal{B}^{\infty}):=(X,\mathcal{B})\times(X,\mathcal{B}%
)\times...$, where $X:=\mathbb{R}^{n}$ and $\mathcal{B}=\mathcal{B}(X)$ is the
Borel $\sigma$-algebra on $X$. Let $\pi(\omega,dx^{\infty})$ be the regular
conditional distribution of $x^{\infty}(\omega)$ given the $\sigma$-algebra
$\mathcal{G}$ (see, e.g., \cite{ArkinEvstigneev1987}, Appendix I, Theorem 1).
Denote by $x_{k}^{\infty}$ the $k$th element of the sequence $x^{\infty
}=(x_{1},x_{2},...)$ regarded as a function of $x^{\infty}$. By virtue of
(\ref{cond-conv}) and in view of the uniform conditional integrability of
$\beta_{k}(\omega):=\rho(x_{k}(\omega),C(\omega))$, we have%
\begin{equation}
\int\pi(\omega,dx^{\infty})x_{k}^{\infty}=E[x_{k}|\mathcal{G}](\omega
)\rightarrow y(\omega)\ \ [x_{k}=x_{k}(\omega)], \label{cond-X-y}%
\end{equation}%
\[
\lim_{H\rightarrow\infty}\sup_{k}\int\pi(\omega,dx^{\infty})\rho(x_{k}%
^{\infty},C(\omega))\mathbf{1}_{\{\rho(x_{k}^{\infty},C(\omega))\geq H\}}%
\]%
\begin{equation}
=\lim_{H\rightarrow\infty}\sup_{k}E[\beta_{k}\mathbf{1}_{\{\beta_{k}\geq
H\}}|\mathcal{G}](\omega)=0 \label{cond-int-X}%
\end{equation}
for all $\omega$ belonging to some $\mathcal{G}$-measurable set $\Omega
_{1}\subseteq\Omega$ of measure one. It follows from (\ref{cond-X-y}) and
(\ref{cond-int-X}) that the assumptions of Theorem 1 are satisfied, and so for
each $\omega\in\Omega_{1}$ there exists a $\mathcal{B}^{\infty}$-measurable
vector function $w^{\omega}(x^{\infty})$ integrable with respect to
$\pi(\omega,\cdot)$ and such that
\[
\int\pi(\omega,dx^{\infty})w^{\omega}(x^{\infty})\leq_{C(\omega)}y(\omega)
\]
and
\[
w^{\omega}(x^{\infty})\in\text{Ls\thinspace}(x^{\infty})\ \text{for\ }%
\pi(\omega,\cdot)\text{-almost all }x^{\infty}%
\]
where Ls\thinspace$(x^{\infty})$ is the set of the limit points of the
sequence $x^{\infty}=(x_{1},x_{2},...)$.

\textit{2nd step.} We will use the following fact. There exists a function
$\psi:[0,1]\times X^{\infty}\rightarrow\mathbb{R}^{n}$\ jointly measurable
with respect to $\mathcal{B}[0,1]\times\mathcal{B}^{\infty}$\ (where
$\mathcal{B}[0,1]$ is the Borel $\sigma$-algebra on $[0,1]$) and possessing
the following property. For each finite measure $\mu$\ on $\mathcal{B}%
^{\infty}$\ and each\textit{\ }$\mathcal{B}^{\infty}$-measurable function
$f:X^{\infty}\rightarrow\mathbb{R}^{n}$,\ there exists $r\in\lbrack0,1]$\ such
that $\psi(r,x^{\infty})=f(x^{\infty})$\ for $\mu$-almost all $x^{\infty}\in
X^{\infty}$\textit{. }This result establishes the existence of a
\textquotedblright universal\textquotedblright\ jointly measurable function
parametrizing all equivalence classes of measurable functions $X^{\infty
}\rightarrow\mathbb{R}^{n}$ with respect to all finite measures: any such
class contains a representative of the form $\psi(r,\cdot)$, where $r$ is some
number in $[0,1]$. The result (extending Natanson \cite{Natanson1961}, Chapter
15, Section 3, Theorem 4) follows from Theorem AI.3 in \cite{EST2004} using
the fact that all uncountable standard measurable spaces are isomorphic to the
segment $[0,1]$ with the Borel $\sigma$-algebra (see e.g. Dynkin and
Yushkevich \cite{DynkinYushkevich1979}, Appendix 2).

\textit{3rd step.} For each $\omega\in\Omega$, consider the set $U(\omega)$ of
those $r\in\lbrack0,1]$ for which the function $\psi(r,\cdot)$ satisfies%
\[
\int\pi(\omega,dx^{\infty})\psi(r,x^{\infty})\leq_{C(\omega)}y(\omega),\;
\]%
\begin{equation}
\psi(r,x^{\infty})\in\text{Ls\thinspace}(x^{\infty})\;\text{for }\pi
(\omega,\cdot)\text{-almost all }x^{\infty}. \label{U}%
\end{equation}
Observe that for $\omega\in\Omega_{1}$\ the set $U(\omega)$ is not empty
because it contains an element $r\in\lbrack0,1]$ such that $\psi(r,x^{\infty
})=w^{\omega}(x^{\infty})$ for $\pi(\omega,\cdot)$-almost all $x^{\infty}$.
Further, the set of pairs $(\omega,r)$ satisfying $r\in U(\omega)$ is
$\mathcal{G}\times\mathcal{B}[0,1]$-measurable because $\pi(\omega,dx^{\infty
})$ is a conditional distribution given $\mathcal{G}$, the function
$\psi(r,x^{\infty})$ is $\mathcal{B}[0,1]\times\mathcal{B}^{\infty}%
$-measurable, $C(\omega)$ and $y(\omega)$ are $\mathcal{G}$-measurable, and
the constraint in (\ref{U}) can be written as
\begin{equation}
\int\pi(\omega,dx^{\infty})F(r,x^{\infty})=1, \label{U1}%
\end{equation}
where $F(r,x^{\infty})$ is the indicator function of the set
\[
\{(r,x^{\infty}):\psi(r,x^{\infty})\in\text{Ls}\,(x^{\infty})\}\in
\mathcal{B}[0,1]\times\mathcal{B}^{\infty}.
\]
The last inclusion follows from the fact that $z\in$Ls$\,(x^{\infty})$ if and
only if for each $M=1,2,...$ and $l=1,2,...$ there exists $k\geq l$ such that
$|z-x_{k}^{\infty}|<1/M$.

\textit{4th step.} By virtue of Aumann's measurable selection theorem (see
above), there exists a $\mathcal{G}$-measurable function $r(\omega)$ such that
$r(\omega)\in U(\omega)$ (a.s.). Define
\[
x(\omega):=\psi(r(\omega),x^{\infty}(\omega))\ \ \ \ [x^{\infty}%
(\omega)=(x_{1}(\omega),x_{2}(\omega),...)].\
\]
Since $\pi(\omega,dx^{\infty})$ is the conditional distribution of $x^{\infty
}(\omega)$ given $\mathcal{G}$ and $r(\omega)$ is $\mathcal{G}$-measurable, we
have%
\[
E[x(\omega)|\mathcal{G}]=E[\psi(r(\omega),x^{\infty}(\omega))|\mathcal{G}%
]=\int\pi(\omega,dx^{\infty})\psi(r(\omega),x^{\infty})\leq y(\omega)\text{
(a.s.)}.
\]
Furthermore, $x(\omega)\in$Ls$\,(x^{\infty}(\omega))$ (a.s.) because this
inclusion is equivalent to the equality $F(r(\omega),x^{\infty}(\omega))=1$
(a.s.) and%
\[
EF(r(\omega),x^{\infty}(\omega))=E\,\{E[F(r(\omega),x^{\infty}(\omega
))|\mathcal{G}]\}=
\]%
\[
E\int\pi(\omega,dx^{\infty})F(r(\omega),x^{\infty})=1
\]
by virtue of (\ref{U}) and (\ref{U1}). Since $x(\omega)\in$Ls$\,(x^{\infty
}(\omega))$ (a.s.), by Aumann's theorem, there exist integer-valued random
variables $1<k_{1}(\omega)<k_{2}(\omega)<...$ such that $\lim x_{k_{m}%
(\omega)}(\omega)=x(\omega)\;$a.s.; this was shown at the end of the proof of
Theorem 1.\hfill$\square\medskip$

\textbf{Acknowledgement.} The authors are grateful to Zvi Artstein, Ilya
Molchanov and Sergey Pirogov for helpful comments.$\medskip$


\begin{thebibliography}{99}                                                                                               %


\bibitem {ArkinEvstigneev1987}Arkin, V.I. and Evstigneev, I.V.,
\textit{Stochastic Models of Control and Economic Dynamics}, Academic Press,
London, 1987.

\bibitem {Artstein1979}Artstein, Z., A note on Fatou's lemma in several
dimensions, \textit{Journal of Mathematical Economics} \textbf{6} (1979) 277--282.

\bibitem {Aumann1966}Aumann, R.J., Existence of competitive equilibria in
markets with a continuum of traders, \textit{Econometrica }\textbf{34} (1966) 1--17.

\bibitem {Balder-Hess1995}Balder, E.J., Hess, C., Fatou's lemma for
multifonctions with unbounded values, \textit{Mathematics of Operations
Research} \textbf{20} (1995) 175--188.

\bibitem {BESH2018}Babaei, E., Evstigneev, I.V. and Schenk-Hopp\'{e}, K.R.,
Von Neumann-Gale dynamics and capital growth in financial markets with
frictions, preprint, 2018.

\bibitem {BabaeiEvstigneevPirogov2018}Babaei, E., Evstigneev, I.V., and
Pirogov, S.A., Stochastic Fixed Points and Nonlinear Perron-Frobenius Theorem,
\textit{Proceedings of the American Mathematical Society} \textbf{146} (2018) 4315-4330.

\bibitem {BahsounEvstigneevTaksar2008}Bahsoun, W., Evstigneev, I.V., and
Taksar, M.I., Rapid paths in von Neumann-Gale dynamical systems,
\textit{Stochastics} \textbf{80} (2008) 129-142.

\bibitem {Cornet-et-al2003}Cornet, B., Topuzu M., and Yildiz, A., Equilibrium
theory with a measure space of possibly satiated consumers, \textit{Journal of
Mathematical Economics} \textbf{39} (2003) 175--196.

\bibitem {DempsterEvstigneevTaksar2006}Dempster, M.A.H., Evstigneev, I.V. and
Taksar, M.I., Asset pricing and hedging in financial markets with transaction
costs: An approach based on the von Neumann-Gale model, \textit{Annals of
Finance} \textbf{2} (2006) 327--355.

\bibitem {DynkinYushkevich1979}Dynkin, E.B., and Yushkevich, A.A.,
\textit{Controlled Markov processes and their applications}, N.Y., Springer, 1979.

\bibitem {EvstigneevSchenkHoppe2008}Evstigneev, I.V., and Schenk-Hopp\'{e},
K.R., Stochastic equilibria in von Neumann-Gale dynamical systems,
\textit{Transactions of the American Mathematical Society} \textbf{360} (2008) 3345--3364.

\bibitem {EST2004}Evstigneev, I.V., Sch\"{u}rger, K., and Taksar, M.I., On the
Fundamental Theorem of Asset Pricing: Random constraints and bang-bang
no-arbitrage criteria, \textit{Mathematical Finance} \textbf{14} (2004) 201-221.

\bibitem {EvstigneevZhitlukhin2013}Evstigneev, I.V., and Zhitlukhin, M.V.,
Controlled random fields, von Neumann-Gale dynamics and multimarket hedging
with risk, \textit{Stochastics} \textbf{85} (2013) 652-666.

\bibitem {Gale1956}Gale, D., A closed linear model of production. In: Kuhn,
H.W. \textit{et al.} (Eds.), \textit{Linear Inequalities and Related Systems},
Princeton University Press, Princeton (1956) 285--303.

\bibitem {Hess2002}Hess, C., Set-valued integration and set-valued probability
theory: an overview. In: Pap, E. (Ed.), \textit{Handbook of Measure Theory},
Elsevier, North--Holland (2002) 617--673.

\bibitem {Hildenbrand1974}Hildenbrand, W., \textit{Core and Equilibria of a
Large Economy}. Princeton University Press, New Jersey, 1974.

\bibitem {Natanson1961}Natanson, I.P., \textit{Theory of Functions of a Real
Variable}, N.Y., Ungar, 1961.

\bibitem {Neveu1965}Neveu, J., \textit{Mathematical Foundations of the
Calculus of Probability}, San Francisco, Holden-Day, 1965.

\bibitem {SaadouneValadier1995}Saadoune, M. and Valadier, M., Extraction of a
\textquotedblright good\textquotedblright\ subsequence from a bounded sequence
of integrable functions. \textit{Journal of Convex Analysis} \textbf{2} (1995) 345-357.

\bibitem {Schmeidler1970}Schmeidler, D., Fatou's lemma in several dimensions,
\textit{Proceedings of the American Mathematical Society} \textbf{24} (1970) 300--306.

\bibitem {vonNeumann1937}Von Neumann, J., \"{U}ber ein \"{o}konomisches
Gleichungssystem und eine Verallgemeinerung des Brouwerschen Fixpunktsatzes,
in:\ \textit{Ergebnisse eines Mathematischen Kolloquiums}, \textbf{8} (1937),
1935--1936 (Franz-Deuticke, Leipzig and Wien), 73--83. [Translated: A model of
general economic equilibrium,\ \textit{Review of Economic Studies} \textbf{13}
(1945-1946) 1--9.]
\end{thebibliography}
\end{document}